\documentclass[english]{article}

\usepackage{amssymb}
\usepackage[cp1251]{inputenc}
\usepackage[english]{babel}
\usepackage{amsmath}
\usepackage{color}
\usepackage{graphics}
\usepackage{epsfig}
\usepackage[all,knot]{xy}
\xyoption{arc}

\def\mapr#1{\smash{\mathop{\buildrel{#1}\over\longrightarrow}}}

\def\Aut{{\hbox{\bf Aut}\;}}
\def\Inn{{\hbox{\bf Inn}\;}}
\def\Der{{\hbox{\bf Der}\;}}
\def\Out{{\hbox{\bf Out}\;}}
\def\id{{\hbox{\bf Id}}}

\def\R{{\bf R}}
\def\cA{\cal A}
\def\cB{\cal B}
\def\cC{{\cal C}}
\def\cD{{\cal D}}
\def\cL{{\cal L}}
\def\cT{{\cal T}}
\def\rg{\frak g}

\newtheorem{theorem}{Theorem}
\newtheorem{definition}{Definition}

\title{Description of coupling in the category of transitive Lie algebroids}

\author{Li XiaoYu \\(Harbin Institute of Technology)\\ A.S.Mishchenko \\(Harbin Institute of Technology, \\Moscow Lomonosov State University)\\
}
\date{10.10.2013}

\begin{document}
\maketitle

\begin{abstract}
In our previous paper (\cite{Mi-Li-2013}) we have given a sufficient and necessary condition when 
the coupling between Lie algebra bundle (LAB) and the tangent bundle exists in the sense of Mackenzie (\cite{Mck-2005}, Definition 7.2.2) for the theory of transitive Lie algebroids.
Namely we have defined a new topology on the group $\Aut(\rg)$ of all automorphisms of the Lie algebra $\rg$, say $\Aut(\rg)^{\delta}$,  and show that tangent bundle $TM$ can be coupled with the Lie algebra bundle $L$  if and only if the Lie algebra bundle L admits a local trivial structure with structural group endowed with such new topology. 

But the question how many couplings exist under these conditions still remains 
open. Here we make the result more accurate and prove that there is a one-to-one correspondence between the family  
$Coup(L)$ 
of all coupling of the Lie algebra bundle $L$
with fixed finite dimensional Lie algebra $\rg$ as the fiber and the structural group $\Aut(\rg)$ of all automorphisms of Lie algebra $\rg$ and the tangent bundle $TM$ and the family $LAB^{\delta}(L)$ of equivalent classes of local trivial structures with structural group
$\Aut(\rg)$ endowed with new topology $\Aut(\rg)^{\delta}$.

This result gives a way for geometric construction of the classifying space for
transitive Lie algebroids with fixed structural finite dimensiaonal Lie algebra
$\rg$. Hence we can clarify a categorical description of the characteristic classes for transitive Lie algebroids and a comparison with that by J. Kubarski (\cite{Kub-91e},\cite{Mi-Li-2012})

\end{abstract}

\section{Introduction}

Transitive Lie algebroids have specific properties that allow to look at the transitive
Lie algebroid as an element of the object of a homotopy functor. Roughly speaking each
transitive Lie algebroids can be described as a vector bundle over the tangent bundle of
the manifold which is endowed with additional structures. Therefore transitive Lie
algebroids admits a construction of inverse image generated by a smooth mapping of
smooth manifolds.

Due to K.Mackenzie
(\cite{Mck-2005})
the construction can be managed as a homotopy functor ${TLA_{\rg}}$ from category of
smooth manifolds to the transitive Lie algebroids. The functor ${TLA_{\rg}}$ associates
with each smooth manifold ${M}$ the set ${TLA_{\rg}(M)}$ of all transitive algebroids with
fixed structural finite dimensional Lie algebra ${\rg}$.

Hence one can construct a classifying space ${\cB_{\rg}}$
(\cite{Mi-2010}\cite{Mi-2011})
such that
the family ${TLA_{\rg}(M)}$ of all transitive  Lie algebroids with fixed Lie algebra ${\rg}$ over the manifold ${M}$
has one-to-one correspondence with the family ${[M,\cB_{\rg}]}$ of homotopy classes of continuous maps :
${
TLA_{\rg}(M)\approx [M,\cB_{\rg}].
}$

In spite of the evident categorical point of view we faced the challenge of geometrical construction of the classifying space, in particular generalization of the Eilenberg-MacLane spaces, realization of the cohomological obstructions for
equivariant mapping and others.

Given a transitive Lie algebroid ${\cA}$ over a manifold ${M}$ the homotopy classification consists of a Lie algebra bundle ${L}$ (LAB) and a linear connection ${\nabla}$ on ${L}$, that is covariant differentiation, that satisfies Leibnitz condition with respect to the fiberwise brackets. The connection ${\nabla}$ should
satisfies the property, that the curvature tensor ${R^{\nabla}}$ is trivial modulo
adjoint operator, that is ${\nabla(mod\, ad)}$ is a {\it coupling} between LAB ${L}$ and tangent bundle ${TM}$.

Therefore for homotopy classification of transitive Lie algebroids there is an open crucial problem: under what condition the coupling
${\nabla(mod\,ad)}$ exists. In  (\cite{Mi-Li-2013}) we represented
a way how to solve the problem of existing the coupling  ${\nabla(mod\, ad)}$ in the terms of the structural group
of the Lie algebra bundles ${L}$.
Namely we have defined a new topology on the group $\Aut(\rg)$ of all automorphisms of Lie algebra $\rg$, say $\Aut(\rg)^{\delta}$,  and show that tangent bundle $TM$ can be coupled with the Lie algebra bundle $L$  if and only if the Lie algebra bundle L admits a local trivial structure with structural group endowed with such new topology.

But the question how many couplings exist under these conditions still remains 
open. Here we make the result more accurate and prove that there is a one-to-one correspondence between the family  
$Coup(L)$ 
of all coupling of the Lie algebra bundle $L$
with fixed finite dimensional Lie algebra $\rg$ as the fiber and the structural group $\Aut(\rg)$ of all automorphisms of Lie algebra $\rg$ and the tangent bundle $TM$ and the family $LAB^{\delta}(L)$ of equivalent classes of local trivial structures with structural group
$\Aut(\rg)$ endowed with new topology $\Aut(\rg)^{\delta}$.

This result gives a way for geometric construction of the classifying space for
transitive Lie algebroids with fixed structural finite dimensiaonal Lie algebra
$\rg$.

As a consequence we can clarify a categorical description of the characteristic classes for transitive Lie algebroids and
 a comparison with that by J. Kubarski (\cite{Kub-91e},\cite{Mi-Li-2012})

\section{Definitions}

Given smooth manifold $M$ consider a vector bundle 
over $TM$ with fiber $\rg$:

$$E\mapr{a}TM\mapr{p_{T}}M.$$

The fiber $\rg$ has the structure of a finite dimensional Lie algebra
and the structural group of the bundle $ E$ is  $ Aut(\rg)$, the group of all automorphisms of the Lie algebra $\rg$. Let
$ p_{E}=p_{T}\cdot a$.
So we have a commutative diagram of two vector bundles
$$
\xymatrix{
E\ar[r]^{a}\ar[d]_{p_{E}}&TM\ar[d]^{p_{T}}\\
M\ar[r]&M
}
$$

The diagram is endowed with additional structure (commutator braces) and then is called
(Mackenzie, definition 3.3.1, Kubarski, definition 1.1.1)
transitive Lie algebroid
$$
\cA = \left\{
\makebox(70,25)[t]{\xymatrix{
E\ar[r]^{a}\ar[d]_{p_{E}}&TM\ar[d]^{p_{T}}\\
M\ar[r]&M
}}; \{\bullet,\bullet\}
\right\}.
$$

The braces $\{\bullet,\bullet\}$ satisfy the natural properties, such that the space $\Gamma^{\infty}(E)$
with braces $\{\bullet,\bullet\}$ forms an infinite dimensional Lie algebra with structure of the $ C^{\infty}(M)$ -- module,  that is

\begin{enumerate}
\item  Skew commutativity: for two smooth sections
    $\sigma_{1},\sigma_{2}\in \Gamma^{\infty}(E)$ one has
    \begin{equation}\label{1}
    \{\sigma_{1},\sigma_{2}\}=-\{\sigma_{2},\sigma_{1}\}\in\Gamma^{\infty}(E),
    \end{equation}
\item Jacobi identity: for three smooth sections
    $\sigma_{1},\sigma_{2},\sigma_{3}\in \Gamma^{\infty}(E)$ one has
    \begin{equation}\label{2}
    \{\sigma_{1},\{\sigma_{2},\sigma_{3}\}\}+
    \{\sigma_{3},\{\sigma_{1},\sigma_{2}\}\}+
    \{\sigma_{2},\{\sigma_{3},\sigma_{1}\}\}=0,
    \end{equation}
\item Differentiation: for two smooth sections $\sigma_{1},\sigma_{2}\in \Gamma^{\infty}(E)$ and smooth function $ f\in C^{\infty}(M)$ one has
    \begin{equation}\label{3}
    \{\sigma_{1},f\cdot\sigma_{2}\}=a(\sigma_{1})(f)\cdot\sigma_{2}+
    f\cdot\{\sigma_{1},\sigma_{2}\}
    \in\Gamma^{\infty}(E).
    \end{equation}

\end{enumerate}

\section{Pullback}

Let $ f:M'\mapr{}M$ be a smooth map. Then one can define an inverse image
(pullback) of the Lie algebroid (Mackenzie, page 156, Kubarski, definition 1.1.4),
$ f^{!!}(\cA)$. Namely, we have the commutative diagram
$$
\xymatrix{
(Tf)^{*}(E)\ar@{=}[r]&E'\ar[r]^{f^{!}}\ar[d]^{a'}&E\ar[d]^{a}\\
&TM'\ar[r]^{Tf}\ar[d]^{p_{T'}}&TM\ar[d]^{p_{T}}\\
&M'\ar[r]^{f}&M
}
$$

Each section $\sigma':M'\mapr{}E$ is induced by the pair $(x',\tau')$
$$
\xymatrix{
&E'\ar[r]^{f^{!}}\ar[d]^{a'}&E\ar[d]^{a}\\
&TM'\ar[r]^{Tf}\ar[d]_{p_{T'}}&TM\ar[d]^{p_{T}}\\
&M'\ar[r]^{f}\ar@/^2pc/[uu]^{\sigma'}\ar@/_/[ruu]_{\tau'}\ar@/_/[u]_{x'}&M\ar@/_2pc/[uu]_{\tau}
}
$$
such that the diagram is commutative. If $\tau'=\tau\circ f$ then two vector fields $ x'$ and $ a\circ\tau$ are related or are interwined by $ f$.
Hence the Lie algebroid structure on $ E'$ is induced by the Lie algebroid structure on $ E$.

This means that given a finite dimensional Lie algebra $\rg$
there is a functor $\cT\cL\cA_{\rg}$ such that
with any manifold $ M$ it assigns the family $\cT\cL\cA_{\rg}(M)$ of all transitive  Lie algebroids
with fixed Lie algebra $\rg$.

The following statement can be proved, (see for example
\cite{Wa-2007})
\begin{theorem}
Each transitive Lie algebroid is locally trivial.
\end{theorem}

This means that for a small neighborhood there is a trivia\-lization
of the vector bundles $ E$, $ TM$, $\ker a=L\approx \rg\times M$ such that
$$
E\approx TM\oplus L,
$$
and the Lie braces are defined by the formula:
$$
[(X,u),(Y,v)]= ([X,Y], [u,v]+X(v)-Y(u)).
$$

\section{Homotopy of pullback}

Using the construction of pullback and the idea by Allen Hatcher
(\cite{Hatcher-2005},Proposition 1.7)
one can prove that the functor $\cT\cL\cA_{\rg}$ is the homotopic functor.
More exactly for two homotopic smooth maps
$ f_{0},f_{1}:M_{1}\mapr{}M_{2}$
and for the transitive  Lie algebroid
$$
\cA=\left(E\mapr{a}TM_{2}\mapr{}M_{2};\{\bullet,\bullet\}\right)
$$
two inverse images $ f_{0}^{!!}(\cA)\, $ and $ f_{1}^{!!}(\cA)\,$ are isomorphic.

\section{Coupling}

Each transitive Lie algebroid

$$
\cA = \left\{
\makebox(70,25)[t]{\xymatrix{
E\ar[r]^{a}\ar[d]_{p^{E}}&TM\ar[d]^{p^{T}}\\
M\ar[r]&M
}}\quad , \{\bullet,\bullet\}
\right\}.
$$
one can represent as an exact sequence of bundles

$$
0\mapr{}L\mapr{}E\mapr{a} TM\mapr{}0.
$$

The bundle $ E $ can be represent as a direct sum of bundles
$$
E=L\oplus TM,
$$

$$
0\mapr{}L\mapr{}
\begin{array}{c}
L \\
\oplus \\
TM
\end{array}
\mapr{} TM\mapr{}0.
$$
\noindent Then each section $ \sigma\in\Gamma(E)$ one can represent as the pair of sections
$$
\sigma = (u,X), \quad u\in\Gamma(L), \quad X\in\Gamma(TM).
$$
Hence the commutator brace for the pair of the sections
$ \sigma_{1}=(u_{1}, X_{1}), \quad \sigma_{2}=(u_{2}, X_{2})$
can be written by the formula

$$
\begin{array}{ll}
\{\sigma_{1},\sigma_{2}\}=&\{(u_{1}, X_{1}),(u_{2}, X_{2})\}=\\\\
&\hskip -1cm=\left(
[u_{1},u_{2}]+\nabla_{X_{1}}(u_{2})-\nabla_{X_{2}}(u_{1})+\Omega(X_{1},X_{2}),
[X_{1},X_{2}]
\right).
\end{array}
$$

Here
$$ \nabla_{X}:\Gamma(L)\mapr{}\Gamma(L)$$
is the covariant gradient of fiberwise differentiation of sections,
where
$$\Omega(X_{1},X_{2})\in \Gamma(L)$$
is classical two-dimensional differential form with values in the fibers of the bundle  $L$.

The covariant derivative for fiberwise differentiation of the sections or so called linear connection in the bundle
$L$
$$ \nabla_{X}:\Gamma(L)\mapr{}\Gamma(L)$$
is an operator that satisfies the following natural conditions:

\begin{enumerate}
\item Fiberwise differentiation with respect to multiplication in the Lie algebra structure of the fibre:
$$
\nabla_{X}([u_{1},u_{2}])=[\nabla_{X}(u_{1}),u_{2}]+
[u_{1},\nabla_{X}(u_{2})],$$
$ u_{1},u_{2}\in\Gamma(L)$.

\item Differentiation of sections in the space  $\Gamma(L)$ as module over the function algebra  $ \cC^{\infty}(M):$
$$
\nabla_{X}(f\cdot u)=X(f)\cdot u+f\cdot\nabla_{X}(u),
$$
$ u\in\Gamma(L)$, $ f\in\cC^{\infty}(M)$.

\item Linear dependence on vector fields:
$$
\nabla_{f\cdot X_{1}+g\cdot X_{2}}=f\cdot \nabla_{X_{1}}+g\cdot \nabla_{X_{2}},
$$
that is
$$
\nabla_{f\cdot X_{1}+g\cdot X_{2}}(u)=
f\cdot \nabla_{X_{1}}(u)+g\cdot \nabla_{X_{2}}(u),
$$
$ X_{1},X_{2}\in\Gamma(TM)$, $ f,g\in\cC^{\infty}(M)$,
$ u\in\Gamma(L)$.
\end{enumerate}

From abstract point of view a covariant derivative can be considered as a pair
$\nabla_{X}=(\cD,X)$,
$$
\cD:\Gamma(L)\mapr{}\Gamma(L), \quad
X:\cC^{\infty}(M)\mapr{}\cC^{\infty}(M),
$$
that satisfies the conditions
\begin{enumerate}
\item $$
\cD([u_{1},u_{2}])=[\cD(u_{1}),u_{2}]+
[u_{1},cD(u_{2})],\quad u_{1},u_{2}\in\Gamma(L).$$

\item
$$
\cD(f\cdot u)=X(f)\cdot u+f\cdot\cD(u),\quad u\in\Gamma(L),
\quad f\in\cC^{\infty}(M).
$$
\end{enumerate}

The association  $\nabla_{X}=(\cD,X)=(\cD(X),X)\,$ satisfies the last condition
\begin{enumerate}
\addtocounter{enumi}{2}
\item $$
\nabla_{f\cdot X_{1}+g\cdot X_{2}}=f\cdot \nabla_{X_{1}}+g\cdot \nabla_{X_{2}},
$$
\end{enumerate}

The family
of all covariant derivatives of fiberwise differentiation is
the space of sections of a transitive Lie algebroid, namely
$\cD_{der}(L)\mapr{}TM$, that is there is a bundle
$\cD_{der}(L)\mapr{}TM$, such that one has the exact sequence

$$
0\mapr{}\Der(L)\mapr{}\cD_{der}(L)\mapr{}TM\mapr{}0.
$$
The bundle $\cD_{der}(L)$ can be constructed as a union of fibres where each fiber
$\cD_{der}(L)_{x}$ in the point $ x\in M$ consists of all covariant derivatives $(\cD,X)$ in the point $ x\in M$.

The exact sequence can be included in the exact diagram

\begin{equation}\label{3}
\xymatrix{
&ZL\ar[r]^{=}\ar[d]^{i}&ZL\ar[d]^{i}\\
&L\ar[r]^{=}\ar[d]^{ad} &L\ar[d]^{ad}\\
0\ar[r]&Aut(L)\ar[r]^{j}\ar[d]^{\natural^{0}}&\cD_{Der}(L)\ar[r]^{a}\ar[d]^{\natural}&TM\ar[r]\ar[d]^{=}&0\\
0\ar[r]&Out(L)\ar[r]^{\bar j}\ar[d]&Out\cD_{Der}(L)\ar[r]^{\bar a}\ar[d]&TM\ar[r]\ar[d]&0\\
&0&0&0
}
\end{equation}

The bundle $E$ can be represent as a direct sum of bundles
$$
E\approx L\oplus TM,
$$
using a splitting $\lambda:TM\mapr{}E$, $a\cdot\lambda=\id$,
$$
\xymatrix{
0\ar[r]&L\ar[r]&E\ar[r]^{a}\ar[d]^{\approx}&TM\ar@/_1.pc/[l]_{\lambda}\ar[r]&0\\
&&L\oplus TM&&
}
$$

Then each section $ \sigma\in\Gamma(E)$ one can represent as the pair of sections
$$
\sigma = (u,X), \quad u\in\Gamma(L), \quad X\in\Gamma(TM).
$$
$$
X=a(\sigma),\quad u=\sigma-\lambda(a(\sigma)).
$$

The commutator brace for the pair of the sections
$\sigma_{1}=(u_{1}, X_{1}), \quad \sigma_{2}=(u_{2}, X_{2})$
can be written by the formula

$$
\begin{array}{ll}
\{\sigma_{1},\sigma_{2}\}=&\{(u_{1}, X_{1}),(u_{2}, X_{2})\}=\\\\
&\hskip -1cm=\left(
[u_{1},u_{2}]+\nabla_{X_{1}}(u_{2})-\nabla_{X_{2}}(u_{1})+\Omega(X_{1},X_{2}),
[X_{1},X_{2}]
\right).
\end{array}
$$

Here
$$
\nabla_{X}:\Gamma(L)\mapr{}\Gamma(L)
$$
is the covariant gradient of fiberwise differentiation of sections,
$$
\nabla_{X}(u)=\{\lambda(X),u\}, \quad u\in\Gamma(L),
$$
$$
\nabla_{X}(f\cdot u)=X(f)\cdot u+f\cdot\nabla_{X}(u),\quad  u\in\Gamma(L),\quad f\in C^{\infty}(M),
$$
$$
\nabla_{X}([u_{1},u_{2}])=[\nabla_{X}(u_{1}),u_{2}]+[u_{1},\nabla_{X}(u_{2})] \quad  u_{1},u_{2}\in\Gamma(L).
$$

The form
$$\Omega(X_{1},X_{2})\in \Gamma(L)$$
is classical two-dimensional differential form with values in the fibers of the bundle  $L$,
$$
\Omega(X_{1},X_{2})=\{\lambda(X_{1}),\lambda(X_{2})\}-\lambda([X_{1},X_{2}]).
$$

The covariant gradient $\nabla$ and the differential form $\Omega$ satisfy the following conditions:

$$
R^{\nabla}(X_{1},X_{2})(u)=-[u,\Omega(X_{1},X_{2})].
$$
and
$$
\begin{array}{l}
d^{\nabla}\Omega(X_{1},X_{2},X_{3})=0.
\end{array}
$$

The splitting $\lambda$ and covariant gradient $\nabla$ are included in the commutative diagram with adjoint homomorphism $ ad$:
$$
\xymatrix{
0\ar[r]&L\ar[r]\ar[d]^{ad}&E\ar[r]\ar[d]^{ad}&TM\ar@/_1.pc/[l]_{\lambda}\ar[r]\ar[d]^{=}&0\\
0\ar[r]&Aut(L)\ar[r]&\cD_{der}(L)\ar[r]&TM\ar@/_1.pc/[l]_{\nabla}\ar[r]&0
}
$$

Using covariant gradient $\nabla$ we extend the exact diagram to
$$
\xymatrix{
&ZL\ar[r]^{=}\ar[d]^{i}&ZL\ar[d]^{i}\\
&L\ar[r]^{=}\ar[d]^{ad} &L\ar[d]^{ad}\\
0\ar[r]&\Der(L)\ar[r]^{j}\ar[d]^{\natural^{0}}&\cD_{Der}(L)\ar[r]^{a}\ar[d]^{\natural}&TM\ar@/_1.pc/[l]_{\nabla}\ar[r]\ar[d]^{=}&0\\
0\ar[r]&\Out(L)\ar[r]^{\bar j}\ar[d]&\cD_{Out}(L)\ar[r]^{\bar a}\ar[d]&TM\ar@/_1.pc/[l]_{\Xi}\ar[r]\ar[d]&0\\
&0&0&0
}
$$

Since the curvature tensor $ R^{\nabla}$ satisfies the condition
$$
R^{\nabla}(X_{1},X_{2})(u)=-[u,\Omega(X_{1},X_{2})]
$$
or
$$
R^{\nabla}=ad\circ\Omega
$$
one has
$$
R^{\Xi}(X_{1},X_{2})\equiv 0.
$$
Such a map $\Xi$ is called a {\it coupling} $ TM$ with $ L$.

\begin{theorem}
The coupling $\Xi$ generated by a transitive Lie algebroid $\cA$ does not depend of the choice of the splitting $\lambda$.
\end{theorem}

\section{Classification}

So the classification of the transitive Lie algebroids can be represented by several steps.

The first step consists of constructing a coupling  between LAB $ L$ and tangent bundle $ TM$, that is a homomorphism $\Xi$
$$
\xymatrix{
0\ar[r]&\Out(L)\ar[r]^{\bar j}&\cD_{Out}(L)\ar[r]^{\bar a}&TM\ar@/_1.pc/[l]_{\Xi}\ar[r]&0
}
$$
that satisfies the condition:
$$
R^{\Xi}(X_{1},X_{2})\equiv 0,
$$
where $ R^{\Xi}$ is the map
$$
R^{\Xi}:\Lambda^{2}(TM)\mapr{}\Out(L),
$$
defined by the formula
$$
R^{\Xi}(X_{1},X_{2})=\{\Xi_{X_{1}},\Xi_{X_{2}}\}-\Xi_{[X_{1},X_{2}]}.
$$

So the first step consists of the problem how to describe all couplings between
LAB $ L$ and tangent bundle $ TM$ in the terms of the homotopy theory.

The all subsequent steps were described by Mackenzie
(\cite{Mck-2005})
as the 3-dimensional
cohomological class of an obstruction for existence of transitive Lie algebroid and 2-dimensional cohomological class as parameter for description of all transitive Lie algebroids with fixed coupling $\Xi$ between
LAB $ L$ and tangent bundle $ TM$.

\section{New topology}

Denote by $\Aut\rg$ the group of Lie algebra automorphisms of $\rg$ and by $\Der\rg$ denotes the Lie algebra of derivations of $\rg$. The subgroup of $\Aut\rg$ generated by $ exp(\psi)$, $\psi\in\Der$ is denoted by $\Inn \rg$ and its elements are called inner automorphisms.

\begin{definition}
Denote by $\Aut^{\delta}\rg$ the space $\Aut\rg$ with finer topology such that topology of $\Aut\rg/\Inn \rg$ becomes discrete topology.
\end{definition}

\section{Existence of coupling}

Let $ L$ be a Lie algebra bundle (LAB) on smooth manifold $ M$ with fibre $ \rg$.
\begin{theorem}

A coupling $$\Xi:TM\rightarrow \Out\cD_{Der}(L)$$ exists if and only if $ L$ admits a locally trivial structure with structural group $ {\Aut^{\delta}\rg}$.
\end{theorem}

\section{Uniqueness of couplings}
\subsection{Classification of couplings}

The natural conjecture consists of that the family of all couplings has a one-to-one correspondence with the family
of all classes of locally trivial structures with structural group $ \Aut^{\delta}\rg.$

Namely, given a LAB
$ L $ denote by $ Coup(L)$  the family of all couplings of the LAB $ L$ with the tangent bundle $ TM$.

Denote by $ LAB^{\delta}(L)$ the family of all classes of locally trivial structures with structural group $ \Aut^{\delta}\rg$ on the LAB $ L$. The conjecture says that there is a one-to-one correspondence
$$
f: Coup(L)\mapr{} LAB^{\delta}(L).
$$

The family $ Coup(L)$ can be described by the classes of linear connections $ \nabla$,
$$
\nabla_{X}: \Gamma^{\infty}(L)\mapr{}\Gamma^{\infty}(L),
$$
that satisfies the conditions
$$
R^{\nabla}(X,Y)(u)\equiv [\Omega(X,Y),u]
$$
for some differential form $\Omega$.
In this case we say that the connection $\nabla$ is in accordance with the structural group
$ \Aut^{\delta}\rg.$
Two linear connections $\nabla$ and $\nabla'$ that are in accordance with the structural group
$ \Aut^{\delta}\rg$ belong to the same class
if there is a map $ l:TM\mapr{} L$ such that
$$
\nabla'_{X}(u)\equiv \nabla_{X}(u)+[l(X),u].
$$

\begin{theorem}
If the connection $\nabla$ is in accordance with the structural group
$ \Aut^{\delta}\rg$ then the linear connection $ \nabla'_{X}(u)\equiv \nabla_{X}(u)+[l(X),u]$
also is in accordance with the structural group
$ \Aut^{\delta}\rg$,
that is satisfies the condition
$$
R^{\nabla'}(X,Y)(u)\equiv [\Omega'(X,Y),u]
$$
for some proper differential form $\Omega'$.
\end{theorem}

\subsection{Classification of $LAB^{\delta}$}

By the definition a locally trivial structures with structural group $ \Aut^{\delta}\rg$ on the LAB $ L$
is given for a sufficiently fine open covering $ \{ U_{\alpha}\}$ as a system of trivializations
$$
\varphi_{\alpha}:U_{\alpha}\times\rg\mapr{}p^{-1}(U_{\alpha})\subset L
$$
that satisfy the conditions: for all $ x\in U_{\alpha\beta}=U_{\alpha}\cap U_{\beta}$
$$
\varphi_{\alpha\beta}(x)=\varphi_{\beta}(x)\varphi_{\alpha}^{-1}(x)\in \Aut(\rg)
$$
and the map
$$
\varphi_{\alpha\beta}: U_{\alpha\beta}\mapr{} \Aut^{\delta}(\rg),
$$
is continuous.

Two locally trivial structures with structural group $ \Aut^{\delta}\rg$ on the LAB $ L$,

$ \varphi_{\alpha}:U_{\alpha}\times\rg\mapr{}p^{-1}(U_{\alpha})\subset L$ and
$ \varphi'_{\alpha}:U_{\alpha}\times\rg\mapr{}p^{-1}(U_{\alpha})\subset L$,

are called equivalent if for all $ x\in U_{\alpha}$
$$
\varphi'^{-1}_{\alpha}(x)\circ\varphi_{\alpha}(x)\in \Aut(\rg)
$$
and the map
$$
\varphi'^{-1}_{\alpha}\circ\varphi_{\alpha}:U_{\alpha}\mapr{} \Aut^{\delta}(\rg)$$
is continuous.

\subsection{Description of the map $f: Coup(L)\mapr{} LAB^{\delta}(L)$}

Let fix an open atlas of charts $\{U_{\alpha}\}$ such that each chart $ U_{\alpha}$ is diffeomorphic to
$ \R^{n}$. Let fix central points $ x_{\alpha}\in U_{\alpha}$ and isomorphisms
$$
\psi_{\alpha}:\rg\mapr{}p^{-1}(x_{\alpha})\subset L.
$$
Consider a system of smooth paths $\gamma_{\alpha,x}$ that start in the point
$ x_{\alpha}$ and terminate in the point $ x\in U_{\alpha}$. As an example of such system is the image of rays in $\R^{n}$ which is diffeomorphic to the chart
$ U_{\alpha}$.

Given a linear connection $\nabla$ consider the map
$$
\varphi^{\nabla}_{\alpha}(x):\rg\mapr{}p^{-1}(x)\subset L, \quad x\in U_{\alpha}
$$
that is the result of parallel transport along the curve $\gamma_{\alpha,x}$ with respect
to the connection $\nabla$.

\begin{theorem}
For $ x\in U_{\alpha\beta}$ the map
$$
\varphi^{\nabla}_{\alpha\beta}(x)=(\varphi^{\nabla}_{\beta})^{-1}(x)\circ
\varphi^{\nabla}_{\alpha}(x):\rg\mapr{}\rg
$$
belongs to $\Aut\rg$ and the map
$$
\varphi^{\nabla}_{\alpha\beta}:U_{\alpha\beta}\mapr{}\Aut^{\delta}\rg
$$
is continuous.
\end{theorem}

So the correspondence $ \nabla \Rrightarrow \{\varphi^{\nabla}_{\alpha}\}$ induce
a map
$$
f: Coup(L)\mapr{} LAB^{\delta}(L), \quad f(\nabla)=\{\varphi^{\nabla}_{\alpha}\}.
$$

\subsection{Inverse of the map $g: LAB^{\delta}(L)\mapr{}Coup(L)$}

The inverse map
$$
g: LAB^{\delta}(L)\mapr{}Coup(L)
$$
is defined by the following:

Given a local trivialization $ \Phi=\{\varphi_{\alpha}\}$ and a partition of units
$ H=\{h_{\alpha}\}$
subordinate to the covering $\{U_{\alpha}\}$
let put
$$
\nabla^{\Phi,H}_{X}(u) = \sum\limits_{\alpha}\nabla^{\Phi,\alpha}_{X}(h_{\alpha}\cdot u),\quad u\in\Gamma^{\infty}(L,M),
$$
where
$$
\nabla^{\Phi,\alpha}_{X}(v(x))=\varphi_{\alpha}(x)\frac{\partial}{\partial X}\varphi_{\alpha}^{-1}(x)(v(x)),
\quad v\in\Gamma^{\infty}(L,U_{\alpha}).
$$

\begin{theorem}
The class $ g(\Phi,H)=\nabla^{\Phi,H}\in Coup(L)$ does not depend of the choice of equivalent local trivialization $\Phi$ and the partition of units $ H$ that is
the formula $ g(\Phi,H)=\nabla^{\Phi,H}$
defines a correct map
$$
g: LAB^{\delta}(L)\mapr{}Coup(L).
$$
\end{theorem}

\begin{theorem}
The map $ g$ is inverse to the map $ f$
$$
\xymatrix{
Coup(L)\ar@/^1.5pc/[r]_{f}&LAB^{\delta}(L)\ar@/^1.5pc/[l]_{g}
}
$$
\end{theorem}

This means that we have a one-to-one correspondence between the family of all couplings
$ Coup(L)$ and the family of equivalent classes of locally trivial structures
on the Lie algebra bundle $ L$ with structural group $\Aut^{\delta}\rg$.

\bibliographystyle{plain}
\bibliography{LiXiaoyu-Mishchenko-CouplingUniqeness-references}

\end{document}